\documentclass[12pt]{article}
\usepackage{bbm}
\usepackage{amsfonts}
\usepackage{pifont}
\usepackage{hyperref}
\usepackage{mathrsfs}
\usepackage{indentfirst}
\usepackage{amsmath}
\usepackage{amssymb}
\usepackage[amsmath, thmmarks]{ntheorem}
\usepackage{cite}
\usepackage{graphicx}
\usepackage{tikz}
\newtheorem{definition}{\bf Definition}[section]
\newtheorem{lemma}{\bf Lemma}[section]
\newtheorem{theorem}{\bf Theorem}[section]
\newtheorem{remark}{\bf Remark}[section]
\newtheorem{corollary}{\bf Corollary}[section]
\newtheorem{proposition}{\bf Proposition}[section]
\newtheorem{example}{\bf Example}[section]

\def\QEDopen{{\setlength{\fboxsep}{0pt}\setlength{\fboxrule}{0.2pt}\fbox{\rule[0pt]{0pt}{1.3ex}\rule[0pt]{1.3ex}{0pt}}}} 
\def\QED{\QEDopen}
\def\proof{{\bf Proof.} }
\def\endproof{\hspace*{\fill}~\QED\par\endtrivlist\unskip}

\begin{document}
\setcounter{page}{1}

\title{{\textbf{Triangular norms on finite lattices}}\thanks {Supported by the National Natural Science
Foundation of China (No. 12471440)}}
\author{Peng He\footnote{\emph{E-mail address}: 443966297@qq.com}, Xue-ping Wang\footnote{Corresponding author. xpwang1@hotmail.com; fax: +86-28-84761502}\\
\emph{\small (School of Mathematical Sciences, Sichuan Normal University}\\ \emph{\small Chengdu 610066, Sichuan, People's Republic of China)}}
\newcommand{\pp}[2]{\frac{\partial #1}{\partial #2}}
\date{}
\maketitle

\begin{quote}
{\bf Abstract}
This article intends to characterize triangular norms on a
finite lattice.
We first give a method for generating a triangular norm on an atomistic lattice by the values
of atoms. Then we prove that every triangular norm on a non-Boolean atomistic lattice is not left-continuous and $T_M$ is the uniquely left-continuous triangular norm on an atomistic Boolean lattice. Furthermore, we show that each atomistic Boolean lattice can be represented by a family of triangular norms on an atomistic lattice with the same number of atoms. Finally, we construct a triangular norm on a finite lattice by restricting a triangular norm on an extended atomistic lattice of the finite lattice to the finite lattice.

{\textbf{\emph{Keywords}}:}\ Triangular norm; Atomistic lattice; Continuity; Boolean lattice; Finite lattice
\end{quote}

\section{Introduction}\label{intro}
In order to generalize the triangle inequality in a metric space to more general probabilistic metric spaces, Menger \cite{Meg} first introduced triangular norms (briefly t-norms) on the closed unit interval $[0,1]$ of the real line. After that, Schweizer and Sklar \cite{sch, Schw} further studied t-norms and presented their axioms that are still used today. Initially, t-norms have been used extensively for defining the intersection of fuzzy sets
and modelling logical `and' in fuzzy logic. In fuzzy logic, the set of truth values is modelled by
the real unit interval. However, in modern fuzzy logic, more generally, the real unit interval is replaced by a bounded lattice, on which there are more practical applications \cite{Kl}. Due to the close connection between fuzzy sets theory and order theory, several authors have worked with t-norms on bounded
lattices \cite{Sun, De, Sam, Coomann94, Pal, Pa14, El, ca}.

Regarding the construction of t-norms on lattices, many scholars mainly use the t-norms of
special sublattices and the transformation of known functions. For example, Saminger \cite{Sam} provided
necessary and sufficient conditions for an ordinal sum operation resulting in again a t-norm on
bounded lattices whereas the operation is determined by an arbitrary selection of
subintervals as carriers for arbitrary summand t-norms.
Saminger et al. \cite{Sam08} also investigated the smallest and largest possible extensions
of t-norms on bounded lattices. Latter, Palmeira and Bedregal \cite{Pal} presented a method of extending t-norms from a sublattice $M$ to a bounded lattice $L$ by considering a more general version of the idea of the sublattice. In particular, Palmeira et al. \cite{Pa14} established a way to extend t-norms to a lattice-valued setting by preserving the largest possible number of properties of those fuzzy connectives which are invariants under homomorphisms. Recently, Sun and Liu \cite{Sun} studied t-norms of bounded lattices by additive generators.

Based on the fact that every element in a finite distributive lattice has a uniquely join-irreducible decomposition, Y{\i}lmaz and Kazanc{\i} \cite{yi} gave a method for constructing t-norms on a finite distributive lattice $L$ by means of join irreducible elements of $L$. The work of Y{\i}lmaz and Kazanc{\i} \cite{yi} motivates us consider whether we can construct t-norms on an atomistic lattice by its atoms or not because all elements of an atomistic lattice can be represented by its atoms. It is well-known that a finite lattice can surely be cover-preserving extended into an atomistic lattice, therefore, is it possible to construct t-norms for all finite lattices? This article will positively answer the above two questions.

The rest of this article is organized as follows. In Section 2, we recall some basic concepts and results required lately. In Section 3, we give a method for generating t-norms on atomistic lattices just by their atoms. In Section 4, we show that every t-norm on a non-Boolean atomistic lattice is not left-continuous and $T_M$ is the uniquely left-continuous t-norm on an atomistic Boolean lattice. In Section 5, we show that every atomistic Boolean lattice can be represented by a family of t-norms on atomistic lattices. In Section 6, we present a method for constructing t-norms on a finite lattice $L$ by means of the restriction of t-norms on an extended atomistic lattice of $L$. A conclusion is drawn in Section 7.

\section{Preliminaries}
In this section, we recall some necessary concepts and facts about lattices and t-norms. For more
comprehensive presentations, see, e.g., \cite{Birkhoff73, Crawley73, Ma, De, sch, Kl}.

A poset is a structure $(P, \leq)$ where $P$ is a nonempty set and $\leq$ is an
ordering (reflexive, antisymmetric and transitive)
relation on $P$. Although the set $P$ is only the domain of the partial order $\leq$, we usually follow the custom of identifying $P$ with the poset. A bounded poset $(P, \leq)$ is a poset that possesses
the bottom element $0$ and top element $1$. A lattice is a poset $(L, \leq)$ in which every two elements
subset $\{x, y\}$ has the greatest lower bound, meet, denoted by $x\wedge y$, and the least upper bound, join, denoted by $x\vee y$. A bounded lattice $(L, \leq)$ is a lattice that possesses
the bottom element $0$ and top element $1$. A poset in which every subset has a join and a meet is called a complete lattice.
Note that each complete lattice is a bounded lattice with $\wedge \emptyset=\vee L=1$ and $\vee \emptyset=\wedge L=0$. An element $c$ in a complete lattice $L$ is called compact if whenever $c\leq \vee S$
there exists a finite subset $T\subseteq S$ with $c\leq \vee T$. We define a lattice $L$ to be algebraic if $L$ is
complete and each of its element is a join of compact elements.

Let $L$ be a lattice. For all $a, b\in L$, $a\|b$ denotes that $a\ngeq b$  and $a \nleq b$, and $a\nparallel b$ denotes that $a\geq b$ or
$a\leq b$. The symbol $a\prec b$ means that $a< b$ and there is no element $c\in L$
such that $a<c<b$. An element that covers the least element $0$ of a lattice $L$ will be
referred to as an atom of $L$. We say that an element $q$ in a lattice $L$ is join irreducible if,
for all $x, y\in L$, $q=x\vee y$ implies $q=x$ or $q=y$. The set of non-zero join-irreducible elements and
the set of atoms of $L$ will be denoted by $J(L)$ and $A(L)$, respectively. Note that $A(L)\subseteq J(L)$ for any lattice $L$. Further, for a fixed lattice $L$, we denote $A(x)=\{a\in A(L)| a\leq x\}$
and $J(x)=\{a\in J(L)| a\leq x\}$ for any $x\in L$.
The length of $L$, i.e., $\mbox{sup}\{n: L\mbox{ has an }(n+1)\mbox{-element chain}\}$,
will be denoted by $\ell(L)$.

Let $A$ and $B$ be two sets. We define $A\setminus B=\{x\in A| x\notin B\}$.
For convenience, if $B=\{b\}$, then we denote $A\setminus B=A\setminus b$.

\begin{definition}[\cite{Birkhoff73}]\label{de000}
\emph{A lattice $L$ is called an atomistic lattice if, for all $x\in L$
there exists a set $S\subseteq A(L)$ such that $x=\vee S$.}
\end{definition}

Obviously, an atomistic lattice $L$ is complete and each atom of $L$ is a compact element.
Thus, every atomistic lattice is an algebraic lattice. Indeed, $\vee A(x)=x$ for any $x\in L$ when
$L$ is atomistic.

\begin{definition}[\cite{Crawley73}]\label{de0001}
\emph{Let $L$ be a complete lattice. A subset $S\subseteq L$ is said to be
independent if, for all $a\in L$, $$a\wedge \vee(S\setminus a)=0.$$}
\end{definition}

If $L$ is an algebraic lattice, then a subset $S$ is independent if and only if every finite subset of $S$
is independent. With this fact in hand, a straightforward application of Zorn's lemma yields the following:
If $A$ is a subset of an algebraic lattice $L$ and $B$ is an independent subset of $A$, then $B$ is contained
in a maximal independent subset of $A$. Specially, if $L$ is an atomistic lattice and $B_{x}$ is a maximal
independent of $A(x)$ then $\vee B_{x}=x$ for any $x\in L$ (cf. \cite{Crawley73}).

\begin{definition}[\cite{Schw}]\label{de002}
\emph{Let $P$ be a bounded poset. A t-norm on $P$ is a binary operation on $P$ that is monotone, commutative, associative and with neutral element $1$.}
\end{definition}

The following are the strongest and the weakest t-norms, respectively, on a bounded lattice $L$:\\
$$T_{M}(x, y)=x\wedge y\mbox{ and } T_{D}(x,y)=\begin{cases}
x\wedge y &\mbox{if } 1\in\{x, y\},\\
0 & \mbox{otherwise}.
\end{cases}$$

\begin{definition}[\cite{De}]\label{de001}
\emph{Suppose that $T_{1}$ and $T_{2}$ are two t-norms on a bounded poset $L$. If $T_{1}(x,y) \leq T_{2}(x,y)$ for all $x, y\in L$, then one writes $T_{1}\leq T_{2}$.}
\end{definition}

Denote by $T(L)$ the set of all t-norms on a bounded poset $L$. Then under the partial order $\leq$ defined as Definition \ref{de001}, $T(L)$ is a poset. In particular, if $L$ is a bounded lattice then $T(L)$ is also a bounded poset, i.e., $T_{D}\leq T\leq T_{M}$ for any $T\in T(L)$ since $0\leq T(x, y)\leq x\wedge y$ for any $x, y\in L$.

\begin{definition}[\cite{De}]\label{de002}
\emph{Let $L$ be a bounded poset and $T\in T(L)$. An element $x\in L$ is called an
idempotent element of $T$ if $T(x,x)=x$. The set of all idempotent elements of $T$ is denoted by $\mathcal{I}_{T}(L)$.}
\end{definition}

\begin{definition}[\cite{De,Ma,Zhang}]\label{de003}
\emph{Let $L$ be a complete lattice and $T\in T(L)$.}\\
\emph{(i) $T$ is said to be left-continuous if for any $a\in L$ and any $S\subseteq L$ the following holds:
$$T(a, \vee S)=\bigvee_{s\in S}T(a,s).$$}\\
\emph{(ii) $T$ is said to be right-continuous if for any $a\in L$ and any $S\subseteq L$ the following holds:
$$T(a, \wedge S)=\bigwedge_{s\in S}T(a,s).$$}\\
\emph{(iii) $T$ is said to be continuous if and only if it is both left-continuous and right-continuous.}
\end{definition}

\section{T-norms generated by atoms of atomistic lattices}
Let $L$ be an atomistic lattice and $C(L)=A(L)\cup \{0, 1\}$. This section develops a method to construct all t-norms on $C(L)$, which are then used to generate t-norms of the atomistic lattice $L$.

It is easy to see that a t-norm on a lattice of $\ell(L)= 1$ is unique which equals $T_M$ and $T_D$. Therefore, in the following we always suppose that $L$ is a lattice of $\ell(L)\geq 2$.

For any bounded lattice $L$, it is clear that $C(L)$ is a sub-poset of $L$. In particular, $C(L)$ is an atomistic lattice
and $1\leq \ell(C(L))\leq 2$. Further, if $\ell(C(L))=2$ then $A(L)\neq \emptyset$. Note that the lattice $C(L)$ may not be a sublattice of $L$, generally. Define $$\kappa(x)=\{u\in C(L)\setminus0 | u\leq x\}$$ for any $x\in L$.
It is easy to see that $\kappa(x)=A(x)$ for any $x\in L\setminus 1$. Specially, if $L$
is atomistic with $\ell(L)\geq 2$ then $$A(L)\neq \emptyset, \ell(C(L))=2 \mbox{ and } 1\notin A(L)$$ and, $\kappa(x)=A(x)$ if and only if  $x\in L\setminus 1$.

Let $\mathcal{P}(X)$ be the set of all subsets of a set $X$. Then we have the following theorem.
\begin{theorem}\label{Thm3.1}
Let $L$ be an atomistic lattice. Then $T\in T(C(L))$ if and only if there exists a set $\alpha\in \mathcal{P}(A(L))$ such that $T=T_{\alpha}$ where $T_{\alpha}$ is defined by
\begin{equation}\label{13}
 T_{\alpha}(x,y)=
\begin{cases}
x\wedge y, &\mbox{if } x=1 \mbox{ or } y=1,\\
x,& x=y \in \alpha,\\
0, &  otherwise.
\end{cases}
\end{equation}
\end{theorem}
\proof
Let $\alpha\in \mathcal{P}(A(L))$ and $T=T_{\alpha}$. Obviously, $T_{\alpha}$ is commutative with neutral element $1$.
Next we show that $T_{\alpha}$ is monotone and associative.

(i) Let $x, y, z\in C(L)$ and $y\leq z$. If $x\in \{0,1\}$ or $y=z$ then clearly $T_{\alpha}(x,y)\leq T_{\alpha}(x,z)$. Suppose $x\notin \{0,1\}$ and $y< z$. Note that $x\in A(L)$ and $y\neq 1$. If $y=0$
then $T_{\alpha}(x,y)=0\leq T_{\alpha}(x,z)$. If $y\neq 0$ then $y\in A(L)$ and $z=1$. Thus $T_{\alpha}(x,y)=0$ or $T_{\alpha}(x,y)=x$, which together with $T_{\alpha}(x,z)=x$ implies that $T_{\alpha}(x,y)\leq T_{\alpha}(x,z)$. Hence, $T_{\alpha}$ is monotone.

(ii) Let $x, y, z\in C(L)$. We prove that $T_{\alpha}$ is associative, i.e., $$T_{\alpha}(x,T_{\alpha}(y,z))=T_{\alpha}(T_{\alpha}(x,y), z).$$
If $1\in \{x,y,z\}$ or $0\in \{x,y,z\}$ then the above equation is obvious. Now, suppose
that $x, y, z\in A(L)$. Then the proof of the associativity is made in three cases.

Case 1. If $x, y$ and $z$ are all different from each other then $T_{\alpha}(x,T_{\alpha}(y,z))=0=T_{\alpha}(T_{\alpha}(x,y), z)$.

Case 2. If $x=y=z$ then $T_{\alpha}(x,T_{\alpha}(y,z))=T_{\alpha}(T_{\alpha}(x,y), z)$ since $T_{\alpha}$ is commutative.

Case 3. If $x=y\neq z$ then $T_{\alpha}(x,T_{\alpha}(y,z))=T_{\alpha}(x,0)=0$. On the other hand, $T_{\alpha}(x,y)$
is either equal to $x$ or equal to $0$. This means that $T_{\alpha}(T_{\alpha}(x,y), z)=0$. Thus $T_{\alpha}(x,T_{\alpha}(y,z))=T_{\alpha}(T_{\alpha}(x,y), z)$. The subcases both $x\neq y=z$ and $x=z\neq y$ are analogous.

From Cases 1, 2 and 3, $T_{\alpha}$ is associative. Therefore, $T=T_{\alpha}\in T(C(L))$.

Conversely, suppose that $T\in T(C(L))$. Then from the construction of $C(L)$, we have that
\begin{equation}\label{02}
T(u,u)=u\mbox{ or }T(u,u)=0 \mbox{ for any }u\in A(L)
\end{equation} and
\begin{equation}\label{03}
T(u,v)=0 \mbox{ for any }(u,v)\in A(L)^2\mbox{ with }u\neq v.
\end{equation}
Set $\alpha=\{u\in A(L)|T(u,u)=u\}$. Then $\alpha\in \mathcal{P}(A(L))$ and $T=T_{\alpha}$.
\endproof
\begin{theorem}\label{The1}
Let $L$ be an atomistic lattice and $T\in T(C(L))$. Then the binary operation $\overline{T}$ on $L$
defined by
\begin{equation}\label{01}
\overline{T}(x,y)=\bigvee_{u\in \kappa(x)}\bigvee_{v\in \kappa(y)}T(u,v)
\end{equation}
is a t-norm on $L$.
\end{theorem}
\proof
For the sake of convenience of the proof, we denote
\begin{equation}\label{04}
\theta_{T}(x)=\{u\in A(x)| T(u,u)=u\} \mbox{ and }\beta_{T}(x)=\{u\in A(x)| T(u,u)=0\}
\end{equation}
for any $x\in L$.
Then $\theta_{T}(x)\cup \beta_{T}(x)=A(x)$ and $\theta_{T}(x)\cap \beta_{T}(x)=\emptyset$ by Eq. \eqref{02}. Now, we prove that $\overline{T}$ is a t-norm on $L$.

(i) $\overline{T}(x,1)=\bigvee_{u\in \kappa(x)}\bigvee_{v\in \kappa(1)}T(u,v)=\bigvee_{u\in \kappa(x)}T(u,1)=\bigvee_{u\in \kappa(x)}u=x$ for any $x\in L$ since $T(u,v)\leq T(u,1)$ for any $v\in C(L)$ and $1\in \kappa(1)$.

(ii) Let $x,y,z\in L$ and $y\leq z$. Then it is clear that $\kappa(y)\subseteq \kappa(z)$. Thus $$\overline{T}(x,y)=\bigvee_{u\in \kappa(x)}\bigvee_{v\in \kappa(y)}T(u,v)\leq \bigvee_{u\in \kappa(x)}\bigvee_{v\in \kappa(z)}T(u,v)=\overline{T}(x,z).$$

(iii) It is easy to see that $\overline{T}(x,y)=\overline{T}(y,x)$ for any $x,y\in L$.

(iv) We prove that $\overline{T}(x,\overline{T}(y,z))=\overline{T}(\overline{T}(x,y), z)$ for any $x,y,z\in L$. Note that $\kappa(x)=A(x)$ for any $x\neq 1$. Then from (i), for any $x,y\in L$ we have
\begin{equation}\label{05}
 \overline{T}(x,y)=\bigvee_{u\in \kappa(x)}\bigvee_{v\in \kappa(y)}T(u,v)=
\begin{cases}
y, &\mbox{if } x=1,\\
x, &\mbox{if } y=1,\\
\bigvee_{u\in A(x)}\bigvee_{v\in A(y)}T(u,v),& \mbox{otherwise},
\end{cases}
\end{equation}
which, together with Eqs. (\ref{02}), (\ref{03}) and (\ref{04}), implies
\begin{equation}\label{06}
 \overline{T}(x,y)=
\begin{cases}
y, &\mbox{if } x=1,\\
x, &\mbox{if } y=1,\\
\bigvee_{u\in \theta_{T}(x\wedge y)}u,& \mbox{otherwise}
\end{cases}
\end{equation}since $A(x)\cap A(y)=A(x\wedge y)$ for any $x, y\in L$. If $1\in\{x,y,z\}$, then obviously $\overline{T}(x,\overline{T}(y,z))=\overline{T}(\overline{T}(x,y), z)$. Now, let $1\notin \{x,y,z\}$. Then by (\ref{06})
\begin{equation}\label{eq06}\overline{T}(x,\overline{T}(y,z))=\overline{T}(x,\bigvee_{u\in \theta_{T}(y\wedge z)}u).\end{equation} Set $\omega=A(\bigvee_{u\in \theta_{T}(y\wedge z)}u)$. Because $\bigvee_{u\in\theta_{T}(y\wedge z)}u\leq y\wedge z$, it follows from (\ref{04}) that
\begin{equation}\label{07}
\theta_{T}(y\wedge z)\subseteq \omega\subseteq A(y\wedge z)=\theta_{T}(y\wedge z)\cup\beta_{T}(y\wedge z).
\end{equation}
Then \begin{align*}\overline{T}(x,\bigvee_{u\in \theta_{T}(y\wedge z)}u)&=\bigvee_{u\in A(x)}\bigvee_{v\in\omega}T(u,v)\\
&=[\bigvee_{u\in A(x)}\bigvee_{v\in\omega}T(u,v)]\vee [\bigvee_{u\in A(x)}\bigvee_{v\in\beta_{T}(y\wedge z)\setminus\omega}T(u,v)]
\end{align*}
since $\bigvee_{u\in A(x)}\bigvee_{v\in\beta_{T}(y\wedge z)\setminus\omega}T(u,v)=\bigvee_{u\in A(x)\cap [\beta_{T}(y\wedge z)\setminus\omega]}T(u,u)=0$ by (\ref{03}) and (\ref{04}). Thus
\begin{align*}
\overline{T}(x,\bigvee_{u\in \theta_{T}(y\wedge z)}u)&=\bigvee_{u\in A(x)}\bigvee_{v\in\omega \cup [\beta_{T}(y\wedge z)\setminus\omega]}T(u,v)\\
 &=\bigvee_{u\in A(x)}\bigvee_{v\in A(y\wedge z)}T(u,v) \ \ \ \ \mbox{ by } (\ref{07})\\
 &=\bigvee_{u\in A(x)\cap A(y\wedge z)}T(u,u) \ \ \ \ \mbox{ by } (\ref{03})\\
 &=\bigvee_{u\in A(x\wedge y\wedge z)}T(u,u),
\end{align*}
which together with \eqref{eq06} yields $\overline{T}(x,\overline{T}(y,z))=\bigvee_{u\in A(x\wedge y\wedge z)}T(u,u)$.
Further, from (iii), we have that
$$\overline{T}(\overline{T}(x,y), z)=\overline{T}(z, \overline{T}(x,y))=\bigvee_{u\in A(z\wedge x\wedge y)}T(u,u).$$
Therefore, $\overline{T}(x,\overline{T}(y,z))=\overline{T}(\overline{T}(x,y), z)$.

To sum up, $\overline{T}$ is a t-norm on $L$.
\endproof

Let $L$ be an atomistic lattice. In what follows, the symbol $\overline{T}$ always stands for the t-norm on $L$ generated by a t-norm $T\in T(C(L))$ through (\ref{01}) and set $$\mathfrak{T}(L)=\{ \overline{T} | T\in T(C(L))\}.$$ Denote by $T|_{K}$ the restriction of a binary function $T: P\times P\rightarrow P$ on the subset $K$ of the set $P$. Then we have the following lemma.

\begin{lemma}\label{lem2}
Let $K$ be a bounded sub-poset of a bounded poset $P$ with the top element $1_{p}\in K$ and $T$ be a t-norm on the poset $P$. Then $T|_{K}$ is a t-norm on $K$ if and only if $T$ is closed on $K$, i.e., $T(x, y)\in K$ for any $x,y \in K$.
\end{lemma}
\proof It is trivial that $T$ is closed on $K$ if $T|_{K}$ is a t-norm on $K$. So that we just show that $T|_{K}$ is a t-norm on $K$ when $T$ is closed on $K$. Since $T$ is a t-norm on $P$ and closed on $K$, $T|_{K}$ is associative, commutative and monotone.
The top element $1_{P}$ is a neutral element of $T|_{K}$ because $T(x,1_{P})=x$ and $1_{P}\in K$. So that $T|_{K}$ is a t-norm on $K$.
\endproof
By Lemma \ref{lem2} and Theorem \ref{The1}, the following remark is obviously.
\begin{remark}\label{rem3}
\emph{Let $L$ be an atomistic lattice and $T\in T(C(L))$. Then $\overline{T}|_{C(L)}= T$.}
\end{remark}

\begin{lemma}\label{lem4}
Let $L$ be an atomistic lattice and $T\in T(C(L))$. Then \\
(i) $\overline{T}=T_{M}$ if and only if $T(u,u)=u$ for any $u\in A(L)$;\\
(ii) $\overline{T}=T_{D}$ if and only if $T(u,u)=0$ for any $u\in A(L)$.
\end{lemma}
\proof
(i) If $\overline{T}=T_{M}$ then $T(u,u)=\overline{T}|_{C(L)}(u,u)=\overline{T}(u,u)=u\wedge u=u$ for any $u\in A(L)$ by Remark \ref{rem3}.

Conversely, suppose that $x,y\in L$ and $T(u,u)=u$ for any $u\in A(L)$. Then by \eqref{04}, $\theta_{T}(x)=A(x)$. So, if both $0$ and $1$ do not
belong to $\{x, y\}$, then $$\overline{T}(x,y)=\bigvee_{u\in \theta_{T}(x\wedge y)}u=\bigvee_{u\in A(x\wedge y)}u=x\wedge y$$ by (\ref{06}). If either $0$ or $1$ belongs to $\{x, y\}$, then $\overline{T}(x,y)=x\wedge y$ obviously.
Therefore, $\overline{T}=T_{M}$.

(ii) By Remark \ref{rem3}, $\overline{T}=T_{D}$ implies that $T(u,u)=\overline{T}|_{C(L)}(u,u)=\overline{T}(u,u)=0$ for any $u\in A(L)$.

Conversely, suppose that $T(u,u)=0$ for any $u\in A(L)$ and $x,y\in L\setminus 1$. Then it is clear that $\overline{T}(x,y)=0$, which implies that $\overline{T}=T_{D}$ since $\overline{T}$ is a t-norm on $L$.
\endproof
By Lemma \ref{lem4}, we have the following remark.

\begin{remark}\label{rem+3}
\emph{Let $L$ be an atomistic lattice. Then $(\mathfrak{T}(L), \leq)$ is a bounded poset.}
\end{remark}

\section{Behaviour of t-norms $\overline{T}$ generated by t-norms $T$}
This section shows that every t-norm on a non-Boolean atomistic lattice is not left-continuous and $T_M$ is the uniquely left-continuous t-norm on an atomistic Boolean lattice.

It is well-known that if $L$ is an atomistic Boolean lattice then $L\cong (\mathcal{P}(A(L)), \subseteq)$. So that, in the following, we can use
the standard notation $\mathbf{2}^{A(L)}$ instead of the atomistic Boolean lattice $L$.

\begin{proposition}\label{pro4}
Let $L$ be an atomistic lattice and $T\in T(C(L))$. Then $x\in \mathcal{I}_{\overline{T}}(L)\setminus \{0,1\}$ if and only if there exists a maximal
independent set $B_{x}$ of $A(x)$ such that $T(u,u)=u$ for any $u\in B_{x}$.
\end{proposition}
\proof
Suppose that $x\in \mathcal{I}_{\overline{T}}(L)\setminus \{0,1\}$. Then by (\ref{06}), $$\overline{T}(x, x)=x=\bigvee_{u\in A(x)}T(u,u)=\bigvee_{u\in \theta_{T}(x)}u.$$ So, there exists a maximal
independent set $B_{x}\subseteq \theta_{T}(x)$ with $\vee B_{x}=x$ since $L$ is an atomistic lattice.
Thus $B_{x}\subseteq \theta_{T}(x)\subseteq A(x)$. Therefore, $B_{x}$ is a maximal independent set of $A(x)$ and $T(u,u)=u$ for any $u\in B_{x}$.

Let $x\in L\setminus\{0,1\}$ and $B_{x}$ be a maximal independent set of $A(x)$ with $T(u,u)=u$ for any $u\in B_{x}$. Then $\vee B_{x}=x$. It follows that $$x\geq \overline{T}(x,x)=\bigvee_{u\in A(x)}T(u,u)\geq\bigvee_{u\in B_{x}}T(u,u)=\bigvee_{u\in B_{x}}u=x,$$ and consequently $\overline{T}(x,x)=x\in \mathcal{I}_{\overline{T}}(L)\setminus \{0,1\}$.
\endproof

\begin{definition}\label{def5}
\emph{Let $T$ be a t-norm on a complete lattice $L$. Then $T$ is said to be left-semicontinuous, if for any $a\in L$ and $S\subseteq L$ with $\vee S\neq 1$, the following equation $$T(a, \vee S)=\bigvee_{x\in S}T(a,x)$$ is established.}
\end{definition}

It is evident that every left-continuous t-norm is left-semicontinuous. Note that the t-norm $T_{D}$ on a complete lattice is left-semicontinuous.

\begin{theorem}\label{the6}
Let $L$ be an atomistic lattice.
Then $T_{\diamond}$ is a left-semicontinuous t-norm on $L$ if and only if there
exists a $T\in T(C(L))$ such that $T_{\diamond}=\overline{T}$ where $T$ satisfies that $T(u,u)=0$ for any $u\in A(\vee S)\setminus \bigcup_{x\in S} A(x)$ while $S\subseteq L$ with $\vee S\neq1$.
\end{theorem}
\proof
Suppose that $T_{\diamond}$ is a left-semicontinuous t-norm on $L$. Let $T=T_{\diamond}|_{C(L)}$. From the structure
of $C(L)$, $T_{\diamond}$ is closed on $C(L)$. Thus $T\in T(C(L))$ by Lemma \ref{lem2}. If $x, y\in L\setminus 1$, then
$$T_{\diamond}(x,y)=T_{\diamond}(\vee A(x), \vee A(y))=\bigvee_{u\in A(x)}\bigvee_{v\in A(y)}T_{\diamond}(u,v)=\bigvee_{u\in A(x)}\bigvee_{v\in A(y)}T(u,v)$$ since $L$ is an atomistic lattice and $T=T_{\diamond}|_{C(L)}$. So $T_{\diamond}=\overline{T}$ by (\ref{05}).

For any $S\subseteq L$ with $\vee S\neq1$, let $u\in A(\vee S)\setminus \bigcup_{x\in S} A(x)$ and suppose that $T(u,u)\neq 0$, i.e., $T(u,u)=u$. Then by (\ref{03}),
\begin{equation}\label{08}
T_{\diamond}(u,\vee S)=\overline{T}(u,\vee S)=\bigvee_{v\in A(\vee S)}T(u,v)=T(u,u)=u\neq 0
\end{equation} since $\vee S\neq1$.
However, $u\notin A(x)$ for any $x\in S$, which together with (\ref{03}) yields that $T_{\diamond}(u,x)=\overline{T}(u,x)=\bigvee_{v\in A(x)}T(u,v)=0$ for any $x\in S$. It follows from Eq. (\ref{08}) that $\bigvee_{x\in S}T_{\diamond}(u,x)=0\neq T_{\diamond}(u,\vee S)$, contrary to the fact that $T_{\diamond}$ is left-semicontinuous. Therefore, $T(u,u)=0$ for any $u\in A(\vee S)\setminus \bigcup_{x\in S} A(x)$.

Now, we prove the sufficiency. Let $S\subseteq L$ with $\vee S\neq 1$. If $a=1$ then obviously
$$T_{\diamond}(a,\vee S)=\overline{T}(a,\vee S)=\overline{T}(1,\vee S)=\vee S=\bigvee_{x\in S}\overline{T}(1,x)=\bigvee_{x\in S}\overline{T}(a,x)=\bigvee_{x\in S}T_{\diamond}(a,x).$$  Next, we suppose that $a\neq1$.
Then \begin{align*}T_{\diamond}(a,\vee S)&=\overline{T}(a,\vee S)\\&=\bigvee_{u\in A(a)}\bigvee_{v\in A(\vee S)}T(u,v)\\
&=\bigvee_{u\in A(a)\cap A(\vee S)}T(u,u) \ \ \ \mbox{ by } (\ref{03})\\
&=\bigvee_{u\in A(a\wedge \vee S)}T(u,u).
\end{align*}
Similarly, we have that
\begin{align*}\bigvee_{x\in S}T_{\diamond}(a,x)&=\bigvee_{x\in S}\overline{T}(a,x)\\&
=\bigvee_{x\in S}\bigvee_{u\in A(a)}\bigvee_{v\in A(x)}T(u,v)\\
&=\bigvee_{x\in S}\bigvee_{u\in A(a)\cap A(x)}T(u,u)\ \ \ \mbox { by } (\ref{03})\\
&=\bigvee_{x\in S}\bigvee_{u\in A(a\wedge x)}T(u,u) \\
&=\bigvee_{u\in \bigcup_{x\in  S}A(a\wedge x)}T(u,u).
\end{align*}
Note that $\bigcup_{x\in  S}A(a\wedge x)\subseteq  A(a\wedge \vee S)$. Suppose $u\in  A(a\wedge \vee S) \setminus \bigcup_{x\in  S}A(a\wedge x)$, i.e., $u\in A(a)$ and $u\in A(\vee S)$, but $u\notin A(a\wedge x)$ for any $x\in S$. We claim that $u\notin A(x)$ for any $x\in S$. Otherwise, there exists an $x\in S$ such that $u\in A(a\wedge x)$ since $u\in A(a)$, a contradiction. Hence, $u\notin \bigcup_{x\in S}A(x)$. So $u\in A(\vee S)\setminus \bigcup_{x\in S}A(x)$, which yields that
\begin{equation}\label{09}
T(u, u)=0
\end{equation}
for any $u\in A(a\wedge \vee S) \setminus \bigcup_{x\in  S}A(a\wedge x)$.

To sum up,
\begin{align*}
T_{\diamond}(a,\vee S)&=\bigvee_{u\in A(a\wedge \vee S)}T(u,u)\\&=\bigvee_{u\in [A(a\wedge \vee S)\setminus \bigcup_{x\in  S}A(a\wedge x)]\cup \bigcup_{x\in  S}A(a\wedge x)}T(u,u)\\
&=[\bigvee_{u\in A(a\wedge \vee S)\setminus \bigcup_{x\in  S}A(a\wedge x)}T(u,u)] \vee [\bigvee_{u\in \bigcup_{x\in  S}A(a\wedge x)}T(u,u)]\\
&=\bigvee_{u\in \bigcup_{x\in  S}A(a\wedge x)}T(u,u) \ \ \ \mbox{ by } (\ref{09})\\
&=\bigvee_{x\in S}T_{\diamond}(a,x),
\end{align*}
thus $T_{\diamond}$ is a left-semicontinuous t-norm on $L$.
\endproof

Note that if $L$ is an atomistic Boolean lattice, then $A(\vee S)\setminus \bigcup_{x\in S} A(x)=\emptyset$ for any $S\subseteq L$. So the following corollary is clearly by Theorem \ref{the6}.

\begin{corollary}\label{cor07}
Let $L$ be an atomistic Boolean lattice and $T\in T(L)$. Then $T$ is left-semicontinuous if and
only if $T\in \mathfrak{T}(L)$.
\end{corollary}

\begin{theorem}\label{the8}
A t-norm $T_{\Delta}$ on an atomistic lattice $L$ is left-continuous if and only if $L\cong \mathbf{2}^{A(L)}$ and $T_{\Delta}=T_{M}$.
\end{theorem}
\proof
The sufficiency is obviously. Now suppose that $T_{\Delta}$ is a left-continuous t-norm on $L$, then there exists a t-norm $T\in T(C(L))$ such that $T_{\Delta}=\overline{T}$ by Theorem \ref{the6}. Thus $$u=T_{\Delta}(u,1)=\bigvee_{v\in A(L)}T_{\Delta}(u,v)=\bigvee_{v\in A(L)}\overline{T}(u,v)=\bigvee_{v\in A(L)}T(u,v)=T(u,u)$$ for any $u\in A(L)$ by Eq. (\ref{03}). It follows from Lemma \ref{lem4} that $T_{\Delta}=\overline{T}=T_{M}$.

Further, for any $S\subseteq L$, we know that if $u\in A(\vee S)$ then
$$u=T_{\Delta}(u,u)=T_{\Delta}(u, \vee S)= \bigvee_{x\in S}T_{\Delta}(u, x)=\bigvee_{x\in S}\bigvee_{v\in A(x)}T_{\Delta}(u, v)=\bigvee_{v\in \bigcup_{x\in S}A(x)}T_{\Delta}(u, v)$$ since $T_{\Delta}$ is left-continuous, $L$ is atomistic and $T_{\Delta}=T_{M}$. Thus $u\in \bigcup_{x\in S}A(x)$ for all $u\in A(\vee S)$. This means that
$A(\vee S)\subseteq \bigcup_{x\in S}A(x)$. Consequently,
\begin{equation}\label{10}
A(\vee S)=\bigcup_{x\in S}A(x)
\end{equation}
for any $S\subseteq L$ since $A(\vee S)\supseteq\bigcup_{x\in S}A(x)$.

Let $\varphi: (L, \leq) \rightarrow (\mathcal{P}(A(L)),\subseteq )$ be defined by
$$\varphi(x)=A(x)$$ for all $x\in L$. It is easy to verify that $\varphi$ preserves
the order and it is an injective map. Suppose $U\in \mathcal{P}(A(L))$. Then $\vee U\in L$ and
$A(\vee U)=\bigcup_{u\in U}A(u)=U$ by (\ref{10}). Thus $\varphi(\vee U)=U$ for
any $U\in \mathcal{P}(A(L))$, which means that $\varphi$ is surjective. Meanwhile, $\varphi^{-1}$ preserves
the order since $\varphi^{-1}(U)=\vee U$. Therefore, $\varphi$ is an isomorphic map between lattices $L$ and $\mathcal{P}(A(L))$, i.e., $L\cong \mathbf{2}^{A(L)}$.
\endproof
\begin{remark}\label{rem9}
\emph{From Theorem \ref{the8}, it is clear that there does not exist a left-continuous t-norm on a non-Boolean atomistic lattice.}
\end{remark}

\section{Representation of atomistic Boolean lattices}
Ray \cite{R} proved that an atomistic Boolean lattice $B$ can be
isomorphically represented by a family of t-noms on $B$, under the
partial order of t-norms. In this section, we prove that each atomistic Boolean lattice $B$ can be characterized by a family of t-norms of
an atomistic lattice $L$ with $|A(L)|=|A(B)|$.

In the following, we denote by $T_{\alpha}$ the t-norm defined by \eqref{13}. Then we have a theorem as below.
\begin{theorem}\label{the9}
Let $L$ be an atomistic lattice. Then $\mathfrak{T}(L)\cong \mathbf{2}^{A(L)}$.
\end{theorem}
\proof
We first prove that
\begin{equation}\label{11}
 (\mathfrak{T}(L),\leq) \cong (T(C(L)),\leq).
\end{equation}
Let $\mu: (T(C(L)),\leq ) \rightarrow(\mathfrak{T}(L), \leq)$ be defined by $$\mu(T)=\overline{T}$$ for all $T\in T(C(L))$. It is clear that $\mu$ is well-defined. If $T_1, T_2\in T(C(L))$ with $T_1\neq T_2$, then
obviously $\mu(T_1)=\overline{T_{1}}\neq \overline{T_{2}}=\mu(T_{2})$. Hence $\mu$ is injective. Moreover, it is clear that for any $\overline{T}\in \mathfrak{T}(L)$ there exists a
\begin{equation}\label{12}
 T\in T(C(L))
\end{equation}
with $T=\overline{T}|_{C(L)}$ such that $\mu(T)=\overline{T}$ by Remark \ref{rem3}, i.e., $\mu$ is surjective. Consequently, the map $\mu$ is a bijection. Below, we only need to prove that both $\mu$ and its inverse are order-preserving.

Set $T_1, T_2\in T(C(L))$ and $T_1\leq T_2$, and observe that application of (\ref{01}) yields that
$\mu(T_1)=\overline{T_{1}}\leq \overline{T_{2}}=\mu(T_{2})$. Thus $\mu$ is order-preserving.
Now suppose that $\overline{T_{1}}\leq \overline{T_{2}}$. Then $\mu^{-1}(\overline{T_{1}})=\overline{T_{1}}|_{C(L)}\leq \overline{T_{2}}|_{C(L)}=\mu^{-1}(\overline{T_{2}})$ by (\ref{12}). Thus the inverse of $\mu$ is also order-preserving. Therefore, the formula (\ref{11}) is correct.

Inasmuch as (\ref{11}), it suffices to show that  $T(C(L))\cong \mathbf{2}^{A(L)}$. Let $\nu: (\mathcal{P}(A(L)),\subseteq)\rightarrow (T(C(L)),\leq )$ be defined by $$\nu(\alpha)=T_{\alpha}$$ for all $\alpha \in \mathcal{P}(A(L))$. It is clear that $\nu$ is well-defined. Suppose that $\alpha, \beta\in \mathcal{P}(A(L))$ and $\alpha\neq \beta$. Then $\nu(\alpha)\neq \nu(\beta)$ clearly. Further, if $\alpha\subseteq \beta$ then $\nu(\alpha)\leq \nu(\beta)$ by (\ref{13}). Thus $\nu$ is injective and order-preserving.

From Formulas (\ref{13}), (\ref{02}), (\ref{03}) and (\ref{04}), we know that if $T\in T(C(L))$ then $$\alpha=\theta_{T}(1)=\{u\in A(L)| T(u,u)=u\}\in \mathcal{P}(A(L))$$ and $T=T_{\alpha}$.
Hence $\nu$ is surjective and \begin{equation}\label{14}\nu^{-1}(T)=\theta_{T}(1)\end{equation} for any $T\in T(C(L))$.

Suppose $T_{1}, T_{2}\in T(C(L))$ and $T_{1}\leq T_{2}$. Then $u\in A(L)$ and $T_{1}(u,u)=u$ imply that $T_{2}(u,u)=u$ because $u=T_{1}(u,u)\leq T_{2}(u,u)\leq u$. So $\theta_{T_{1}}(1)\subseteq \theta_{T_{2}}(1)$. It follows from (\ref{14}) that $\nu^{-1}(T_{1})\subseteq \nu^{-1}(T_{2})$. Therefore, $\nu^{-1}$ is order-preserving.

To sum up, $\nu$ is an isomorphic map between $\mathcal{P}(A(L))$ and $T(C(L))$. Thus $T(C(L))\cong \mathbf{2}^{A(L)}$, which together with (\ref{11}) yields that $\mathfrak{T}(L)\cong \mathbf{2}^{A(L)}$.
\endproof

\begin{remark}\label{rem10}
\emph{Let $\alpha, \beta\subseteq A(L)$ with $\alpha \cap \beta=\emptyset$ and $\alpha \cup \beta=A(L)$. Then
$T_{\alpha}, T_{\beta}\in T(C(L))$ and $\overline{T_{\alpha}}, \overline{T_{\beta}}\in \mathfrak{T}(L)$ by Theorems \ref{Thm3.1} and \ref{The1}. Additionally, both the complement elements of $T_{\beta}$ and of $\overline{T_{\beta}}$ are $T_{\alpha}$ and $\overline{T_{\alpha}}$ in the lattices $T(C(L))$ and $\mathfrak{T}(L)$, respectively. Furthermore, $\overline{T_{\alpha}}\vee \overline{T_{\beta}}=\overline{T_{\alpha\cup \beta}}$ and $\overline{T_{\alpha}}\wedge \overline{T_{\beta}}=\overline{T_{\alpha\cap \beta}}$.}
\end{remark}

From Theorem \ref{the9}, we have the following corollary.
\begin{corollary}\label{cor10}
Let $L$ be an atomistic Boolean lattice. Then $\mathfrak{T}(L)\cong L$.
\end{corollary}

\section{T-norms on finite lattices}

This section establishes a construction method of t-norm on a finite lattice $L$ by restricting a t-norm on an extended atomistic lattice $\tilde{L}$ of $L$ to the finite lattice $L$.

Let $L$ be a finite lattice and $H(L)=J(L)\setminus A(L)$. For each $p\in H(L)$, insert a new element $w_{p}$ into $L$ such that $0\prec w_{p}\prec p$ and $w_{p}\neq w_{q}$ while $p\neq q$ where $q\in H(L)$. Let $\tilde{L}=L\bigcup \{w_{p}| p\in H(L)\}$. Then $A(\tilde{L})=\{w_{p}| p\in H(L)\}\cup A(L)$ and the construction of $\tilde{L}$ is depicted in Fig.1 where the black-filled elements are the new elements.

\par\noindent\vskip50pt
\begin{minipage}{11pc}
\setlength{\unitlength}{0.75pt}\begin{picture}(600,180)
\put(140,80){\circle{6}}\put(136,68){\makebox(0,0)[l]
{\footnotesize $0$}}
\put(140,120){\circle{6}}\put(145,118){\makebox(0,0)[l]
{\footnotesize $b$}}
\put(140,200){\circle{6}}\put(138,210){\makebox(0,0)[l]
{\footnotesize $1$}}
\put(180,160){\circle{6}}\put(185,158){\makebox(0,0)[l]
{\footnotesize $c$}}
\put(100,160){\circle{6}}\put(85,158){\makebox(0,0)[l]
{\footnotesize $d$}}
\put(140,83){\line(0,1){34}}
\put(142,122){\line(1,1){36}}
\put(138,122){\line(-1,1){36}}
\put(102,162){\line(1,1){36}}
\put(178,162){\line(-1,1){36}}
\put(130,45){$L$}
\put(340,80){\circle{6}}\put(336,68){\makebox(0,0)[l]
{\footnotesize $0$}}
\put(340,120){\circle{6}}\put(325,118){\makebox(0,0)[l]
{\footnotesize $b$}}
\put(340,200){\circle{6}}\put(338,210){\makebox(0,0)[l]
{\footnotesize $1$}}
\put(380,160){\circle{6}}\put(385,158){\makebox(0,0)[l]
{\footnotesize $c$}}
\put(300,160){\circle{6}}\put(285,158){\makebox(0,0)[l]
{\footnotesize $d$}}
\put(300,120){\circle{6}}\put(280,118){\makebox(0,0)[l]
{\footnotesize $w_{d}$}}
\put(296,116){$\bullet$}
\put(380,120){\circle{6}}\put(385,118){\makebox(0,0)[l]
{\footnotesize $w_{c}$}}
\put(376,116){$\bullet$}
\put(340,83){\line(0,1){34}}
\put(300,123){\line(0,1){34}}
\put(342,122){\line(1,1){36}}
\put(338,122){\line(-1,1){36}}
\put(338,122){\line(-1,1){36}}
\put(302,162){\line(1,1){36}}
\put(378,162){\line(-1,1){36}}
\put(338,82){\line(-1,1){36}}
\put(342,82){\line(1,1){36}}
\put(380,123){\line(0,1){35}}
\put(333,45){$\tilde{L}$}
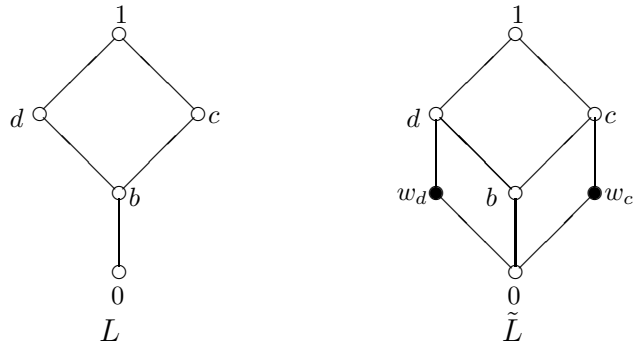
\put(110,20){ Fig.1. The lattices $L$ and $\tilde{L}$, respectively}
\end{picture}
\end{minipage}

\begin{theorem}[\cite{He}]\label{the10}
Let $L$ be a finite lattice. Then $\tilde{L}$ is a finite atomistic lattice and $L$ is a cover-preserving sublattice of $\tilde{L}$ and $\{0_{\tilde{L}}, 1_{\tilde{L}}\}\subseteq L$.
\end{theorem}

From the above theorem, we know that $0_{L}=0_{\tilde{L}}$  and $1_{L}=1_{\tilde{L}}$. So that, for convenience, we always use
$0$ and $1$ instead of $0_{\tilde{L}}$ and $1_{\tilde{L}}$, respectively. Further, since $L$ is a cover-preserving sublattice of $\tilde{L}$, we always use $\wedge$ and $\vee$ instead of $\wedge_{\tilde{L}}$ and $\vee_{\tilde{L}}$, respectively.

\begin{theorem}\label{the11}
Let $L$ be a finite lattice and $\overline{T}\in \mathfrak{T}(\tilde{L})$. Then $\overline{T}|_{L}\in T(L)$ if and only if $T$ satisfies the condition ($\mathcal{C}$).

($\mathcal{C}$) For any $p\in H(L)\setminus 1$, if $T(w_{p},w_{p})\neq 0$ then there exists a $u\in A_{\tilde{L}}(p)\setminus w_{p}$ such that $T(u,u)\neq 0$.
\end{theorem}
\proof
Let $\overline{T}|_{L}\in T(L)$. Then $\overline{T}$ is closed on $L$ by Lemma \ref{lem2}. Thus, if $p\in H(L)\setminus 1$, then
\begin{align*}
\overline{T}|_{L}(p,p)&=\overline{T}(p,p)\\&=\bigvee_{u\in A_{\tilde{L}}(p)}\bigvee_{v\in A_{\tilde{L}}(p)}T(u,v)\\&=\bigvee_{u\in A_{\tilde{L}}(p)}T(u,u) \ \ \ \mbox{ by } \eqref{03}\\
&=T(w_{p}, w_{p})\vee \bigvee_{u\in A_{\tilde{L}}(p)\setminus w_{p}}T(u,u)\\&\neq w_{p}
\end{align*}
since $w_{p}\notin L$. As $w_{p}\in A(\tilde{L})$ and $T\in T(C(\tilde{L}))$, if $T(w_{p},w_{p})\neq 0$ then $T(w_{p},w_{p})= w_{p}$. This follows from $T(w_{p}, w_{p})\vee \bigvee_{u\in A_{\tilde{L}}(p)\setminus w_{p}}T(u,u)\neq w_{p}$ that $\bigvee_{u\in A_{\tilde{L}}(p)\setminus w_{p}}T(u,u)\neq 0$.
So, there exists a $u\in A_{\tilde{L}}(p)\setminus w_{p}$ such that $T(u,u)\neq 0$.

Conversely, suppose that $\overline{T}|_{L}$ is not a t-norm on $L$. Then from Lemma \ref{lem2} and Theorem \ref{the10}, $\overline{T}$ is not closed on $L$.
Then there exist two elements $x, y\in L$ such that
$\overline{T}(x,y)\notin L$. Obviously, $x\neq 1$, otherwise
$\overline{T}(x,y)=y\in L$, a contradiction. Similarly, $y\neq 1$. Thus, there
exists a $p\in H(L)\setminus 1$ such that $$x\wedge y\geq \overline{T}(x,y)=w_{p}=\bigvee_{u\in A_{\tilde{L}}(x)}
\bigvee_{v\in A_{\tilde{L}}(y)}T(u,v)=\bigvee_{u\in A_{\tilde{L}}(x)\cap A_{\tilde{L}}(y)}T(u,u)$$
since $\overline{T}\in \mathfrak{T}(\tilde{L})$. So, $w_{p}\in  A_{\tilde{L}}(x)\cap A_{\tilde{L}}(y)$ and
\begin{equation}\label{15}
T(w_{p}, w_{p})=w_{p}.
\end{equation}
It is clear that $x\wedge y\geq p\succ w_{p}$ in $\tilde{L}$ since $x, y\in L$, $x\wedge y\geq w_{p}$ and $p\in H(L)\setminus 1$. Then $A_{\tilde{L}}(p)\subseteq A_{\tilde{L}}(x\wedge y)=A_{\tilde{L}}(x)\cap A_{\tilde{L}}(y)$. Note that $\overline{T}(w_{p}, w_{p})=T(w_{p}, w_{p})$.
Hence,$$w_{p}=\overline{T}(w_{p}, w_{p})\leq \overline{T}(p,p)
=T(w_{p}, w_{p})\vee \bigvee_{u\in A_{\tilde{L}}(p)\setminus w_{p}}T(u,u) \leq \overline{T}(x,y)=w_{p}.$$
Therefore, $\overline{T}(w_{p}, w_{p})=T(w_{p}, w_{p})\vee \bigvee_{u\in A_{\tilde{L}}(p)\setminus w_{p}}T(u,u)=w_{p}$,
and it follows from \eqref{15} that $\bigvee_{u\in A_{\tilde{L}}(p)\setminus w_{p}}T(u,u)=0$ since $w_{p}$ is an atom of $\tilde{L}$. Consequently, $T(u,u)=0$ for any $u\in A_{\tilde{L}}(p)\setminus w_{p}$, which together with (\ref{15}) implies that the condition ($\mathcal{C}$) does not hold, a contradiction.
\endproof

The following remark is obvious.
\begin{remark}
\emph{Let $L$ be a finite lattice. Then $$\{T\in T(C(\tilde{L}))|T \mbox{ satisfies the condition }(\mathcal{C})\}\neq\emptyset.$$}
\end{remark}
\begin{example}
\emph{Consider the finite lattice $L$ and its extended atomistic lattice $\tilde{L}$ as shown in Fig.1. From Theorem \ref{Thm3.1}, we first construct a t-norm $T$ on $C(\tilde{L})$ as shown by Table \ref{Tab:01}. Then by Theorem \ref{The1}, we obtain a t-norm $\overline{T}$ on $\tilde{L}$ as shown by Table \ref{Tab:02}. Finally, since $T$ satisfies the condition ($\mathcal{C}$), by Theorem \ref{the11} we get a t-norm $\overline{T}|_{L}$ on $L$ as shown by Table \ref{Tab:03}.}
\begin{table}[htpp]
\centering
\caption{The t-norm $T$ on $C(\tilde{L})$}
\label{Tab:01}
\begin{tabular}{c|ccccc}

 $T$ & $0$ & $b$ & $w_d$ & $w_c$ & $1$ \\
 \hline
$0$ & $0$ & $0$ & $0$ & $0$ & $0$ \\
$b$ & $0$ & $b$ & $0$ & $0$ & $b$ \\
$w_d$ & $0$ & $0$ & $w_d$ & $0$ & $w_d$ \\
$w_c$ & $0$ & $0$ & $0$ & $0$ & $w_c$ \\
$1$ & $0$ & $b$ & $w_d$ & $w_c$ & $1$ \\
\end{tabular}
\end{table}

\begin{table}[htpp]
\centering
\caption{The t-norm $\overline{T}$ on $\tilde{L}$}
\label{Tab:02}
\begin{tabular}{c|ccccccc}

 $\overline{T}$ & $0$ & $b$ & $w_d$ & $w_c$ & $d$ & $c$ & $1$ \\
 \hline
$0$ & $0$ & $0$ & $0$ & $0$ & $0$ & $0$ & $0$ \\
$b$ & $0$ & $b$ & $0$ & $0$ & $b$ & $b$ & $b$ \\
$w_d$ & $0$ & $0$ & $w_d$ & $0$ & $w_d$ & $0$ & $w_d$ \\
$w_c$ & $0$ & $0$ & $0$ & $0$ & $0$ & $0$ & $w_c$ \\
$d$ & $0$ & $b$ & $w_d$ & $0$ & $d$ & $b$ & $d$ \\
$c$ & $0$ & $b$ & $0$ & $0$ & $b$ & $b$ & $c$ \\
$1$ & $0$ & $b$ & $w_d$ & $w_c$ & $d$ & $c$ & $1$ \\
\end{tabular}
\end{table}

\begin{table}[htpp]
\centering
\caption{The t-norm $\overline{T}|_{L}$ on $L$}
\label{Tab:03}
\begin{tabular}{c|ccccc}

 $\overline{T}|_{L}$ & $0$ & $b$ & $d$ & $c$ & $1$ \\
 \hline
$0$ & $0$ & $0$ & $0$ & $0$ & $0$ \\
$b$ & $0$ & $b$ & $b$ & $b$ & $b$ \\
$d$ & $0$ & $b$ & $d$ & $b$ & $d$ \\
$c$ & $0$ & $b$ & $b$ & $b$ & $c$ \\
$1$ & $0$ & $b$ & $d$ & $c$ & $1$ \\
\end{tabular}
\end{table}
\end{example}

\begin{corollary}\label{ary}
Let $L$ be a finite lattice with $1\in J(L)$, $\alpha\subseteq \beta\subseteq A(\tilde{L})$ with $\beta\setminus\alpha=\{w_{1}\}$. If $T_{\alpha}$ satisfies the condition ($\mathcal{C}$), then $\overline{T_{\alpha}}|_{L}=\overline{T_{\beta}}|_{L}$.
\end{corollary}
\proof
From the construction of $\tilde{L}$, we know that $0\prec w_{1}\prec 1$ and $w_{1}\notin L$.
Then $T_{\beta}$  satisfies the condition ($\mathcal{C}$) since $\alpha\subseteq \beta\subseteq A(\tilde{L})$, $\beta\setminus\alpha=\{w_{1}\}$ and $T_{\alpha}$
satisfies the condition ($\mathcal{C}$). Thus by Theorem \ref{the11}, both $\overline{T_{\alpha}}|_{L}$ and $\overline{T_{\beta}}|_{L}$ are t-norms on $L$. Now, we prove that
$$\overline{T_{\alpha}}|_{L}(x,y)=\overline{T_{\beta}}|_{L}(x,y)$$ for any $x, y\in L$.
Indeed, if $1\in \{x,y\}$ then $\overline{T_{\alpha}}|_{L} (x,y)=x\wedge y=\overline{T_{\beta}}|_{L}(x,y)$ obviously. Set $x,y\in L\setminus 1$. Thus,
\begin{align*}
\overline{T_{\alpha}}|_{L}(x,y)&=\overline{T_{\alpha}}(x,y)\\&=
\bigvee_{u\in A_{\tilde{L}}(x)}\bigvee_{v\in A_{\tilde{L}}(y)}T_{\alpha}(u,v)\ \ \mbox{ by } (\ref{05})\\
&=\bigvee_{u\in A_{\tilde{L}}(x)\cap A_{\tilde{L}}(y)}T_{\alpha}(u,u)\\
&=\bigvee_{u\in A_{\tilde{L}}(x)\cap A_{\tilde{L}}(y)\cap \alpha}u, \ \ \ \ \mbox{ by } (\ref{13}).
\end{align*}
Similarly, $$\overline{T_{\beta}}|_{L}(x,y)=\bigvee_{u\in A_{\tilde{L}}(x)\cap A_{\tilde{L}}(y)\cap \beta}u.$$
Note that $w_{1}\notin  A_{\tilde{L}}(x)\cap A_{\tilde{L}}(y)$ since $0\prec w_{1}\prec 1$ and $x,y\in L\setminus 1$. Hence, $A_{\tilde{L}}(x)\cap A_{\tilde{L}}(y)\cap \alpha=A_{\tilde{L}}(x)\cap A_{\tilde{L}}(y)\cap \beta$, which means that  $\overline{T_{\alpha}}|_{L}(x,y)=\overline{T_{\beta}}|_{L}(x,y)$.

To sum up, $\overline{T_{\alpha}}|_{L}=\overline{T_{\beta}}|_{L}$.
\endproof
\begin{theorem}\label{the12}
Let $L$ be a finite lattice and $\overline{T}$ be a left-semicontinuous t-norm on $\tilde{L}$. If $\overline{T}|_{L}$ is a t-norm on $L$ then it is left-semicontinuous.
\end{theorem}
\proof
By Lemma \ref{lem2}, $\overline{T}$ is closed on $L$ since $\overline{T}|_{L}$ is a t-norm on $L$.
Thus $\overline{T}|_{L}(x, y)=\overline{T}(x,y)$ for any $x, y\in L$.
Suppose that $x, y, z\in L$ and $y\vee z\neq 1$. Then
$\overline{T}|_{L}(x, y\vee z)=\overline{T}(x, y\vee z)=\overline{T}(x, y)\vee \overline{T}(x, z)=\overline{T}|_{L}(x, y)\vee \overline{T}|_{L}(x, z)$.
Thus $\overline{T}|_{L}$ is left-semicontinuous.
\endproof

\begin{theorem}\label{the13}
Let $L$ be a finite lattice and $\overline{T}$ be a left-semicontinuous t-norm on $\tilde{L}$. If $1\in J(L)$ and $\overline{T}|_{L}$ is a t-norm on $L$ then $\overline{T}|_{L}$ is left-continuous.
\end{theorem}
\proof
By Theorem \ref{the12}, $\overline{T}|_{L}$ is a left-semicontinuous on $L$.
Now, set $x, y, z\in L$ and $y\vee z=1$. Then $y=1$ or $z=1$ since $1\in J(L)$, say $y=1$.
Note that $\overline{T}|_{L}(x,z)\leq x\wedge z\leq x$. Thus, $\overline{T}|_{L}(x, y\vee z)=\overline{T}|_{L}(x, 1)=x=\overline{T}|_{L}(x, 1)\vee \overline{T}|_{L}(x, z)$.
Therefore, $\overline{T}|_{L}$ is left-continuous.
\endproof

\begin{remark}\label{cor14}
\emph{Let $L$ be a finite lattice and $1\in J(L)$. Then there exists at least one left-continuous t-norm on $L$.
Indeed, the t-norm $T_{D}$ is left-continuous.}
\end{remark}
\begin{theorem}\label{theorem}
Let $L$ be a finite lattice and $$\mathcal{S}=\{\overline{T}|_{L}|\overline{T}\in \mathfrak{T}(\tilde{L})\mbox{ and } T \mbox{ satisfies the condition }(\mathcal{C})\}.$$ Then $(\mathcal{S}, \leq)$ is a lattice and $\mathcal{S}\subseteq T(L)$.
\end{theorem}
\proof
From Theorem \ref{the11}, it is obvious that $\mathcal{S}\subseteq T(L)$. From Theorem \ref{Thm3.1},
$\overline{T}|_{L}\in \mathcal{S}$ if
and only if there exists an $\alpha\subseteq A(\tilde{L})$ such
that $T_{\alpha}=T$ on $C(\tilde{L})$ and $T_{\alpha}$ satisfies the condition $(\mathcal{C})$.
Note that $T_{A(\tilde{L})}, T_{\emptyset}\in T(C(\tilde{L}))$ and both of them satisfy the condition $(\mathcal{C})$.
Then $\overline{T_{A(\tilde{L})}}|_{L}, \overline{T_{\emptyset}}|_{L}\in \mathcal{S}$. Furthermore, from Lemma \ref{lem4} and Theorem \ref{the10}, $$\overline{T_{A(\tilde{L})}}|_{L}(x,y)=\overline{T_{A(\tilde{L})}}(x,y)=x\wedge y$$
for any $x,y\in L$, and $$\overline{T_{\emptyset}}|_{L}(x,y)=\overline{T_{\emptyset}}(x,y)=0$$ for
any $x,y\in L\setminus 1$.
Therefore,  $(\mathcal{S}, \leq)$ is a finite bounded poset
with the top element $\overline{T_{A(\tilde{L})}}|_{L}$ and the bottom element $\overline{T_{\emptyset}}|_{L}$.

Let $\overline{T_{\alpha}}|_{L}, \overline{T_{\beta}}|_{L}\in \mathcal{S}$.
Then we know that $\overline{T_{\alpha\cup\beta}}|_{L}\in \mathcal{S}$ since $T_{\alpha\cup\beta}$ satisfies the condition $(\mathcal{C})$ and $\overline{T_{\alpha\cup\beta}}\in \mathfrak{T}(\tilde{L})$. Obviously,
$\overline{T_{\alpha}}|_{L}\leq \overline{T_{\alpha\cup\beta}}|_{L}$ and $\overline{T_{\beta}}|_{L}\leq \overline{T_{\alpha\cup\beta}}|_{L}$.
Next, we prove that $$\overline{T_{\alpha}}|_{L}\vee \overline{T_{\beta}}|_{L}=\overline{T_{\alpha\cup\beta}}|_{L}.$$

Suppose that there exists a $\overline{T_{\gamma}}|_{L}\in \mathcal{S}$ such that $\overline{T_{\alpha}}|_{L}\leq\overline{T_{\gamma}}|_{L}$, $\overline{T_{\beta}}|_{L} \leq\overline{T_{\gamma}}|_{L}$ and $\overline{T_{\gamma}}|_{L}\ngeq\overline{T_{\alpha\cup\beta}}|_{L}$.
Then we claim that $\gamma\nsupseteq \alpha\cup\beta$. Otherwise, $\gamma\supseteq \alpha\cup\beta$ yields that
$\overline{T_{\gamma}}|_{L}\geq\overline{T_{\alpha\cup\beta}}|_{L}$, a contradiction.
We further claim that there exists a $u\in (\alpha\cup\beta)\setminus \gamma$ such that $u\neq w_{1}$ where $0\prec w_{1}\prec 1$.
In fact, if $u=w_{1}$ implies $1\in J(L)$, and it follows from Corollary \ref{ary} that $\overline{T_{\gamma}}|_{L}
=\overline{T_{\gamma\cup\{w_{1}\}}}|_{L}$. Then, begin the whole argument with $\gamma\cup\{w_{1}\}$ instead of $\gamma$ and, there always exists a $u\in (\alpha\cup\beta)\setminus \gamma$ such that $u\neq w_{1}$.
Thus $u\in \alpha$ or $u\in \beta$, say $u\in \alpha$. There are two cases as follows.

Case 1. If $u\in L$ then $\overline{T_{\alpha}}|_{L}(u,u)=\overline{T_{\alpha}}(u,u)=T_{\alpha}(u,u)=u>0=T_{\gamma}(u,u)=\overline{T_{\gamma}}(u,u)
=\overline{T_{\gamma}}|_{L}(u,u)$ since $u\notin \gamma$, contrary to the fact that
$\overline{T_{\alpha}}|_{L}\leq\overline{T_{\gamma}}|_{L}$.

Case 2. If $u\notin \tilde{L}\setminus L$ then there exists a $p\in H(L)\setminus 1$ such that $u=w_{p}$.
Thus there exists a $v\in A_{\tilde{L}}(p)\setminus w_{p}$ such that $v\in \alpha$ since $T_{\alpha}$
satisfies the condition $(\mathcal{C})$. Note that $p=w_{p} \vee v$ in $\tilde{L}$ since $w_{p}\prec p$
and $w_{p}\parallel v$. As $u=w_{p}$ and $v\in \alpha$, we have that
$$p\geq \overline{T_{\alpha}}|_{L}(p,p)=\overline{T_{\alpha}}(p,p)\geq T_{\alpha}(u,u)\vee T_{\alpha}(v,v)=T_{\alpha}(w_{p},w_{p})\vee T_{\alpha}(v,v)=w_{p}\vee v=p,$$
which means that $\overline{T_{\alpha}}|_{L}(p,p)=p$ since $u\notin \gamma$. However, $T_{\gamma}(u,u)=T_{\gamma}(w_{p},w_{p})=0$,
which implies that
\begin{align*}
\overline{T_{\gamma}}|_{L}(p,p)&=\overline{T_{\gamma}}(p,p)\\&=\bigvee_{u\in A_{\tilde{L}}(p)}\bigvee_{v\in A_{\tilde{L}}(p)}T_{\gamma}(u,v)\\
&=\bigvee_{u\in A_{\tilde{L}}(p)}T_{\gamma}(u,u)\\&=\bigvee_{u\in A_{\tilde{L}}(p)\setminus w_{p}}T_{\gamma}(u,u)\\
&\leq \bigvee_{u\in A_{\tilde{L}}(p)\setminus w_{p}} u.
\end{align*}
Thus from the construction of $\tilde{L}$, $\overline{T_{\gamma}}|_{L}(p,p)\leq\bigvee_{u\in A_{\tilde{L}}(p)\setminus w_{p}} u< p=\overline{T_{\alpha}}|_{L}(p,p)$, contrary to the fact that
$\overline{T_{\alpha}}|_{L}\leq\overline{T_{\gamma}}|_{L}$.

To sum up, if $\overline{T_{\alpha}}|_{L}\leq\overline{T_{\gamma}}|_{L}$ and $\overline{T_{\beta}}|_{L} \leq\overline{T_{\gamma}}|_{L}$ then $\overline{T_{\gamma}}|_{L}\geq \overline{T_{\alpha\cup\beta}}|_{L}$. Therefore, $\overline{T_{\alpha}}|_{L}\vee \overline{T_{\beta}}|_{L}=\overline{T_{\alpha\cup\beta}}|_{L}$,
which, together with the poset $\mathcal{S}$ being bounded, yields that $(\mathcal{S},\leq)$ is a lattice.
\endproof

\begin{remark}\label{mark}
\emph{$\mathcal{S}$ is a $\vee$-sublattice of $\mathfrak{T}(\tilde{L})$, but it may not be a sublattice of $\mathfrak{T}(\tilde{L})$ since $\overline{T_{\alpha\cap\beta}}|_{L}\notin \mathcal{S}$, generally.}
\end{remark}

\section{Conclusions}
In this article, we characterized t-norms on a lattice $L$ through the t-norms on a fixed sub-poset of $L$. Since the idea is to follow from part-to-whole, we especially considered atomistic lattices. We proposed a construction method of t-norm on a finite lattice $L$ by restricting a t-norm on an extended atomistic lattice $\tilde{L}$ of $L$ to the finite lattice $L$. From Theorems \ref{Thm3.1} and \ref{The1}, one can
construct $\mathbf{2}^{|A(L)|}$ t-norms on a given atomistic lattice $L$. This implies that the number of left-semicontinuous t-norms on an atomistic Boolean lattice $L$ is exactly $\mathbf{2}^{|A(L)|}$ by Corollary \ref{cor07}.

\end{document}